\documentclass[a4paper,11pt]{amsart}
\usepackage{graphicx}
\usepackage{mathptmx}
\usepackage{amsmath}
\usepackage{amssymb}
\usepackage{enumitem}
\usepackage{xcolor}

\newmuskip\pFqmuskip

\newcommand*\pFq[6][8]{%
  \begingroup 
  \pFqmuskip=#1mu\relax
  \mathcode`=\string"8000
  \begingroup\lccode`\~=`\,
  \lowercase{\endgroup\let~}\pFqcomma
  F^{#2}_{#3}{\left(\genfrac..{0pt}{}{#4}{#5}\bigg|#6\right)}%
  \endgroup
}
\newcommand{\pFqcomma}{\mskip\pFqmuskip}

\newtheorem{theorem}{Theorem}
\newtheorem{lemma}[theorem]{Lemma}
\newtheorem{corollary}[theorem]{Corollary}

\begin{document}

\title[Some relations of two type 2 polynomials and discrete harmonic numbers]{Some relations of two type 2 polynomials and discrete harmonic numbers and polynomials}

\author{Taekyun  Kim}
\address{Department of Mathematics, Kwangwoon University, Seoul 139-701, Republic of Korea}
\email{tkkim@kw.ac.kr}
\author{Dae San  Kim}
\address{Department of Mathematics, Sogang University, Seoul 121-742, Republic of Korea}
\email{dskim@sogang.ac.kr}


\maketitle

\begin{abstract}

The aim of this paper is twofold. The first one is to find several relations between the type 2 higher-order degenerate Euler polynomials and the type 2 high-order Changhee polynomials in connection with the degenerate Stirling numbers of both kinds and Jindalrae-Stirling numbers of both kinds. The second one is to introduce the discrete harmonic numbers and some related polynomials and numbers, and to derive their explicit expressions and an identity.
\end{abstract}

\section{Introduction}

The degenerate Bernoulli and degenerate Euler polynomials were studied by Carlitz in [1], as degenerate versions of the usual Bernoulli and Euler polynomials with their arithmetic and combinatorial interest. The present authors and their colleagues have drawn their attention to various degenerate versions of quite a few special numbers and polynomials and have discovered many properties of them. To name a few, these include the degenerate Stirling numbers of the first and second kinds, degenerate central factorial numbers of the second kind, degenerate Bernoulli numbers of the second kind, degenerate Bell numbers and polynomials, degenerate central Bell numbers and polynomials, degenerate complete Bell polynomials and numbers, degenerate Cauchy numbers, degenerate Bernstein polynomials and so on (see\,\,[2-6,10-12,14,1517-19,22,23]). They have been explored by means of different methods such as generating functions, umbral calculus, combinatorial methods, differential equations, probability theory, $p$-adic integrals., $p$-adic $q$-integrals and special functions. \\
\indent The aim of this paper is twofold. The first one is to find several relations between the type 2 higher-order degenerate Euler polynomials and the type 2 high-order Changhee polynomials in connection with the degenerate Stirling numbers of both kinds and Jindalrae-Stirling numbers of both kinds. We note here that we use the orthogonality relations of the degenerate Stirling numbers in order to derive several corollaries from the obtained theorems. The second one is to introduce the discrete harmonic numbers and the related polynomials and numbers, namely the higher-order degenerate harmonic polynomials and the generalized degenerate harmonic numbers, and to derive their explicit expressions and an identity. \\
\indent This paper is organized as follows. In Section 1, we will recall the stuffs that are needed throughout paper. These include the discrete exponential functions, the discrete Stirling numbers of both kinds, the discrete logarithm function, the harmonic numbers and their generating function, the degenerate Bell polynomials, Jindalrae-Stirling numbers of both kinds, the type 2 higher-order degenerate Euler polynomials and the type 2 higher-order Changhee polynomials.
In Section 2, we show the orthogonality relations for the degenerate Stirling numbers from which the inversion theorem is derived. Then we prove relations between the type 2 degenerate Euler polynomials of order $r$ and the type 2 Changhee polynomials of order $r$ in connection with the degenerate Stirling numbers of both kinds and the Jindalrae-Stirling numbers of both kinds. In additon, we derive a recurrence relation for the type 2 higher-order degenerate Euler polynomials. In Section 3, we introduce the discrete harmonic numbers. Then, as a natural generalization of these numbers, we introduce the higher-order degenerate harmonic polynomials and an explicit expression for those polynomials. Finally, we introduce the generalized degenerate harmonic numbers and find an identity relating these numbers, the type 2 higher-order Changhee polynomials and the degenerate Stirling numbers of the first kind.

\vspace{0.1cm}

For $0 \neq \lambda\in\mathbb{R}$, the degenerate exponential functions are defined by
\begin{equation}
e^{x}_{\lambda}(t)=(1+\lambda t)^{\frac{x}{\lambda}}=\sum_{n=0}^{\infty}(x)_{n,\lambda}\frac{t^{n}}{n!},\quad e_{\lambda}(t)=e_{\lambda}^{1}(t),\quad(\mathrm{see}\
[11,17],\label{1}
\end{equation}
where $(x)_{0,\lambda}=1$, $(x)_{n,\lambda}=x(x-\lambda)\cdots\big(x-(n-1)\lambda\big)$, $(n\ge 1)$. \\
Note that $\displaystyle\lim_{\lambda\rightarrow 0}e_{\lambda}^{x}(t)=e^{xt}\displaystyle$. It is well known that the Stirling numbers of first kind are defined by
\begin{equation}
    (x)_{n}=\sum_{l=0}^{n}S_{1}(n,l)x^{l},\quad(n\ge 0),\quad(\mathrm{see}\ [1-25]), \label{2}
\end{equation}
where $(x)_{0}=1$, $(x)_{n}=x(x-1)\cdots(x-n+1)$, $(n\ge 1)$. \\
As an inversion formula of \eqref{2}, the Stirling numbers of the second kind are defined by
\begin{equation}
    x^{n}=\sum_{l=0}^{n}S_{2}(n,l)(x)_{l},\quad(n\ge 0),\quad(\mathrm{see}\ [7,9-13,16,17]). \label{3}
\end{equation}
As degenerate versions of \eqref{2} and \eqref{3},  the degenerate Stirling numbers of the first kind are given by
\begin{equation}
(x)_{n}=\sum_{l=0}^{n}S_{1,\lambda}(n,l)(x)_{l,\lambda},\ (n\ge 0),\quad (\mathrm{see}\ [12]), \label{5}
\end{equation}
and the degenerate Stirling numbers of the second kind are defined by
\begin{equation}
(x)_{n,\lambda}=\sum_{l=0}^{n}S_{2,\lambda}(n,l)(x)_{l},\quad(n\ge 0),\quad(\mathrm{see}\ [10]). \label{4}
\end{equation}
Let $\log_{\lambda}(t)$ be the compositional inverse of $e_{\lambda}^{x}(t)$, called the discrete logarithm function, such that $\log_{\lambda}(e_{\lambda}(t))=e_{\lambda}\big(\log_{\lambda}(t)\big)=t$. Then we have
\begin{equation}
\log_{\lambda}(1+t)=\frac{1}{\lambda}((1+t)^{\lambda}-1)=\sum_{n=1}^{\infty}\lambda^{n-1}(1)_{n,1/\lambda}\frac{t^{n}}{n!},\quad(\mathrm{see}\ [12]).\label{6}
\end{equation}
The harmonic numbers are defined as
\begin{equation}
H_{0}=1,\quad H_{n}=1+\frac{1}{2}+\frac{1}{3}+\cdots+\frac{1}{n},\quad(n\ge 1),\quad(\mathrm{see}\ [7,9,25]),\label{7}
\end{equation}
We note that the generating function of the harmonic numbers are given by
\begin{equation}
\frac{-\log(1-t)}{1-t}=\sum_{n=1}^{\infty}H_{n}t^{n},\quad(\mathrm{see}\ [7,9,25]). \label{8}
\end{equation}
In [17], the degenerate Bell polynomials are defined by
\begin{equation} \label{9}
e_{\lambda}^x(e_{\lambda}(t)-1)=\sum_{n=0}^{\infty}B_{n,\lambda}(x)\frac{t^n}{n!}.
\end{equation}
When $x=1,\,B_{n,\lambda}=B_{n,\lambda}(1)$ are called the degenerate Bell numbers.
From \eqref{5} and \eqref{4}, we can derive the generating functions of the degenerate Stirling numbers of both kinds which are given by
\begin{equation}
    \frac{1}{k!}\big(e_{\lambda}(t)-1\big)^{k}=\sum_{n=k}^{\infty}S_{2,\lambda}(n,k)\frac{t^{n}}{n!},\quad(\mathrm{see}\ [10]), \label{10}
\end{equation}
and
\begin{equation}
    \frac{1}{k!}\big(\log_{\lambda}(1+t)\big)^{k}=\sum_{n=k}^{\infty}S_{1,\lambda}(n,k)\frac{t^{n}}{n!},\quad(\mathrm{see}\ [10]), \label{11}
\end{equation}
where $k$ is a nonnegative integer. \\
In [16], the Jindalrae-Stirling numbers of the first kind are defined by
\begin{equation}
    \frac{1}{k!}\big(\log_{\lambda}(1+\log_{\lambda}(1+t))\big)^{k}=\sum_{n=k}^{\infty}S_{J,\lambda}^{(1)}(n,k)\frac{t^{n}}{n!}. \label{12}
\end{equation}
As an inversion formula of \eqref{12}, the Jindalrae-Stirling numbers of the second kind are given by
\begin{equation}
    \frac{1}{k!}\big(e_{\lambda}(e_{\lambda}(t)-1)-1\big)^{k}=\sum_{n=k}^{\infty}S_{J,\lambda}^{(2)}(n,k)\frac{t^{n}}{n!},\quad(\mathrm{see}\ [16]). \label{13}
\end{equation}
Carlitz considered the degenerate Euler polynomials of order $r$ which are given by
\begin{equation}
    \bigg(\frac{2}{e_{\lambda}(t)+1}\bigg)^{r}e_{\lambda}^{x}(t)=\sum_{n=0}^{\infty}E_{n,\lambda}^{(r)}(x)\frac{t^{n}}{n!},\quad(\mathrm{see}\ [1]). \label{14}
\end{equation}
The type 2 degenerate Euler polynomials of order $r$ are defined by
\begin{align}
    \bigg(\frac{2}{e_{\lambda}(t)+e_{\lambda}^{-1}(t)}\bigg)^{r}e_{\lambda}^{x}(t)\ &=\ \underbrace{\mathrm{sech}_{\lambda}(t)\times\cdots\mathrm{sech}_{\lambda}(t)}_{r-\mathrm{times}}e_{\lambda}^{x}(t) \label{15}\\
    &=\ \sum_{n=0}^{\infty}\mathcal{E}_{n,\lambda}^{(r)}(x)\frac{t^{n}}{n!},\quad(r\in\mathbb{N}),\nonumber
\end{align}
where
\begin{displaymath}
    \mathrm{sech}_{\lambda}(t)=\frac{2}{e_{\lambda}(t)+e_{\lambda}^{-1}(t)}=\frac{1}{\cosh_{\lambda}(t)},\quad(\mathrm{see}\ [3]).
\end{displaymath}
In [13], the type 2 Changehee polynomials of order $r$ are defined by
\begin{equation}
    \bigg(\frac{2}{(1+t)+(1+t)^{-1}}\bigg)^{r}(1+t)^{x}=\sum_{n=0}^{\infty}C_{n}^{(r)}(x)\frac{t^{n}}{n!}.\label{16}
\end{equation}
When $x=0$, $C_{n}^{(r)}=C_{n}^{(r)}(0)$ are called the type 2 Changhee numbers of order $r$. \\

\section{Type 2 higher-order degenerate Euler and type 2 high-order Changhee polynomials}
From \eqref{5} and \eqref{4}, we note that
\begin{align}
(x)_{n,\lambda}\ &=\ \sum_{k=0}^{n}S_{2,\lambda}(n,k)(x)_{k} \label{17} \\
&= \sum_{k=0}^{n}S_{2,\lambda}(n,k)\sum_{l=0}^{k}S_{1,\lambda}(k,l)(x)_{l,\lambda} \nonumber \\
&= \sum_{l=0}^{n}\bigg(\sum_{k=l}^{n}S_{2,\lambda}(n,k)S_{1,\lambda}(k,l)\bigg)(x)_{l,\lambda}. \nonumber
\end{align}
Therefore, by comparing the coefficients on both sides of \eqref{17}, we obtain the following orthogonality relations where the second one follows analogously to \eqref{17}.
\begin{lemma}
For any integers $n,l$ with $n \geq l$, we have
\begin{displaymath}
\sum_{k=l}^{n}S_{2,\lambda}(n,k)S_{1,\lambda}(k,l)=\left\{\begin{array}{ccc}
1, & \textrm{if $l=n$,}\\
0, & \textrm{otherwise,}
\end{array}\right.
\end{displaymath}
and
\begin{displaymath}
\sum_{k=l}^{n}S_{1,\lambda}(n,k)S_{2,\lambda}(k,l)=\left\{\begin{array}{ccc}
1, & \textrm{if $l=n$,}\\
0, & \textrm{otherwise.}
\end{array}\right.
\end{displaymath}
\end{lemma}
For $n\ge 0$, let $\displaystyle f_{n,\lambda}=\sum_{k=0}^{n}g_{k,\lambda}S_{1,\lambda}(n,k)\displaystyle$.
Then we have
\begin{align}
    \sum_{k=0}^{n}f_{k,\lambda}S_{2,\lambda}(n,k)\ &=\ \sum_{k=0}^{n}\sum_{l=0}^{k}g_{l,\lambda}S_{1,\lambda}(k,l)S_{2,\lambda}(n,k) \label{18} \\
    &=\ \sum_{l=0}^{n}g_{l,\lambda}\bigg(\sum_{k=l}^{n}S_{2,\lambda}(n,k)S_{1,\lambda}(k,l)\bigg)\ =\ g_{n,\lambda}. \nonumber
\end{align}
Therefore, we obtain the following inversion theorem where the converse follows similarly to \eqref{18}.
\begin{theorem}
    Let $n$ be a nonnegative integer. Then we have
    \begin{displaymath}
        f_{n,\lambda}\ =\ \sum_{k=0}^{n}g_{k,\lambda}S_{1,\lambda}(n,k)\ \Longleftrightarrow\ g_{n,\lambda}=\sum_{k=0}^{n}f_{k,\lambda}S_{2,\lambda}(n,k).
    \end{displaymath}
\end{theorem}
In \eqref{15}, replacing $t$ by $\log_{\lambda}(1+t)$, we get
\begin{align}
    \bigg(\frac{2}{(1+t)+(1+t)^{-1}}\bigg)^{r}(1+t)^{x}\ &=\ \sum_{k=0}^{\infty}\mathcal{E}_{k}^{(r)}(x)\frac{1}{k!}\big(\log_{\lambda}(1+t)\big)^{k} \label{19}\\
    &=\ \sum_{k=0}^{\infty}\mathcal{E}_{k}^{(r)}(x)\sum_{n=k}^{\infty}S_{1,\lambda}(n,k)\frac{t^{n}}{n!}\nonumber \\
    &=\ \sum_{n=0}^{\infty}\bigg(\sum_{k=0}^{n}\mathcal{E}_{k}^{(r)}(x)S_{1,\lambda}(n,k)\bigg)\frac{t^{n}}{n!}\nonumber
\end{align}
Therefore, by \eqref{16} and \eqref{19}, we obtain the following theorem.
\begin{theorem}
For $n\ge 0$, we have
\begin{displaymath}
C_{n}^{(r)}(x)\ =\ \sum_{k=0}^{n}\mathcal{E}_{k,\lambda}^{(r)}(x)S_{1,\lambda}(n,k),\quad(r\in\mathbb{N}).
\end{displaymath}
In particular,
\begin{displaymath}
C_{n}^{(r)}\ = \ \sum_{k=0}^{n}\mathcal{E}_{k,\lambda}^{(r)}S_{1,\lambda}(n,k),\quad(r\in\mathbb{N}).
\end{displaymath}
\end{theorem}
By using Theorem 2, we obtain the following corollary.
\begin{corollary}
For $n\ge 0$, we have
\begin{displaymath}
\mathcal{E}_{n,\lambda}^{(r)}(x)\ =\ \sum_{k=0}^{n}C_{k}^{(r)}(x)S_{2,\lambda}(n,k),\quad(r\in\mathbb{N}).
\end{displaymath}
In particular,
\begin{displaymath}
\mathcal{E}_{n,\lambda}^{(r)}\ =\ \sum_{k=0}^{n}C_{k}^{(r)}S_{2,\lambda}(n,k),\quad(r\in\mathbb{N}).
\end{displaymath}
\end{corollary}
Now, we observe that
\begin{align}
\sum_{n=0}^{\infty}\mathcal{E}_{n,\lambda}^{(r)}(x+2)\frac{t^{n}}{n!}+\sum_{n=0}^{\infty}\mathcal{E}_{n,\lambda}^{(r)}(x)\frac{t^{n}}{n!}\ &=\ \bigg(\frac{2}{e_{\lambda}(t)+e_{\lambda}^{-1}(t)} \bigg)^{r}e_{\lambda}^{x}(t)\big(e_{\lambda}^{2}(t)+1\big) \label{20} \\
&=\ \bigg(\frac{2}{e_{\lambda}(t)+e_{\lambda}^{-1}(t)} \bigg)^{r}e_{\lambda}^{x+1}(t)\big(e_{\lambda}(t)+e_{\lambda}^{-1}(t)\big). \nonumber \\
&=\ 2\bigg(\frac{2}{e_{\lambda}(t)+e_{\lambda}^{-1}(t)}\bigg)^{r-1}e_{\lambda}^{x+1}(t)=2\sum_{n=0}^{\infty}\mathcal{E}_{n,\lambda}^{r-1}(x+1)\frac{t^{n}}{n!}.\nonumber
\end{align}
Thus, we have
\begin{equation}
\mathcal{E}_{n,\lambda}^{(r)}(x+2)+\mathcal{E}_{n,\lambda}^{(r)}(x)=2\mathcal{E}_{n,\lambda}^{(r-1)}(x+1), \quad(n\ge 0,\,r\ge 2).\label{21}
\end{equation}
From \eqref{21}, we note that
\begin{align}
\mathcal{E}_{n,\lambda}^{(r)}(x+2)\ &=\ -   \mathcal{E}_{n,\lambda}^{(r)}(x)+2  \mathcal{E}_{n,\lambda}^{(r-1)}(x+1)\label{22} \\
    &=\ -\mathcal{E}_{n,\lambda}^{(r)}(x)+2\big(-\mathcal{E}_{n,\lambda}^{(r-1)}(x-1)+2\mathcal{E}_{n,\lambda}^{(r-2)}(x)\big) \nonumber \\
    &=\ -\mathcal{E}_{n,\lambda}^{(r)}(x)-2 \mathcal{E}_{n,\lambda}^{(r-1)}(x-1)+2^{2}\mathcal{E}_{n,\lambda}^{(r-2)}(x)\nonumber \\
        &=\ \cdots \nonumber \\
    &=\ -\sum_{l=0}^{r-1}2^{l}\mathcal{E}_{n,\lambda}^{(r-l)}(x-l)+2^{r}(x-r+2)_{n,\lambda}. \nonumber
\end{align}
Therefore, by \eqref{22}, we obtain the following theorem.
\begin{theorem}
    For $n\ge 0,\ r\in\mathbb{N}$, we have
    \begin{displaymath}
        \mathcal{E}_{n,\lambda}^{(r)}(x+2)+\sum_{l=0}^{r-1}2^{l}\mathcal{E}_{n,\lambda}^{(r-l)}(x-l)\ =\ 2^{r}(x-r+2)_{n,\lambda}.
    \end{displaymath}
Replacing $t$ by $\log_{\lambda}\big(1+\log_{\lambda}(1+t)\big)$ in \eqref{15}, we get
\begin{align}
&\bigg(\frac{2}{\big(1+\log_{\lambda}(1+t)\big)+\big(1+\log_{\lambda}(1+t)\big)^{-1}}\bigg)^{r}\big(1+\log_{\lambda}(1+t)\big)^{x}\label{23} \\
&=\ \sum_{k=0}^{\infty}\mathcal{E}_{k,\lambda}^{(r)}(x)\frac{1}{k!}\big(\log_{\lambda}(1+\log_{\lambda}(1+t))\big)^{k} \nonumber\\
&=\sum_{k=0}^{\infty}\mathcal{E}_{k,\lambda}^{(r)}(x)\sum_{n=k}^{\infty}S_{J,\lambda}^{(1)}(n,k)\frac{t^{n}}{n!} \nonumber\\
&=\ \sum_{n=0}^{\infty}\bigg(\sum_{k=0}^{n}\mathcal{E}_{n,\lambda}^{(r)}(x)S_{J,\lambda}^{(1)}(n,k)\bigg)\frac{t^{n}}{n!}.\nonumber
\end{align}
\end{theorem}
On the other hand, by \eqref{16}, we get
\begin{align}
&\bigg(\frac{2}{\big(1+\log_{\lambda}(1+t)\big)+\big(1+\log_{\lambda}(1+t)\big)^{-1}}\bigg)^{r}\big(1+\log_{\lambda}(1+t)\big)^{x}  \label{24} \\
&=\ \sum_{k=0}^{\infty}C_{k}^{(r)}(x)\frac{1}{k!}\big(\log_{\lambda}(1+t)\big)^{k}\nonumber \\
&=\ \sum_{k=0}^{\infty}C_{k}^{(r)}(x)\sum_{n=k}^{\infty}S_{1,\lambda}(n,k)\frac{t^{n}}{n!}\nonumber \\
&=\ \sum_{n=0}^{\infty}\bigg(\sum_{k=0}^{n}C_{k}^{(r)}(x)S_{1,\lambda}(n,k)\bigg)\frac{t^{n}}{n!}\nonumber
\end{align}
Therefore, by \eqref{23} and \eqref{24}, we obtain the following theorem.
\begin{theorem}
    For $n\ge 0$, we have
    \begin{displaymath}
        \sum_{k=0}^{n}C_{k}^{(r)}(x)S_{1,\lambda}(n,k)\ =\ \sum_{k=0}^{n}\mathcal{E}_{k,\lambda}^{(r)}(x)S_{J,\lambda}^{(1)}(n,k).
    \end{displaymath}
In particular,
\begin{displaymath}
    \sum_{k=0}^{n}C_{k}^{(r)}S_{1,\lambda}(n,k)\ =\ \sum_{k=0}^{n}\mathcal{E}_{k,\lambda}^{(r)}S_{J,\lambda}^{(1)}(n,k).
\end{displaymath}
\end{theorem}
From Theorem 2, we obtain the following corollary.
\begin{corollary}
    For $n\ge 0$, we have
    \begin{displaymath}
        C_{n}^{(r)}(x)\ =\ \sum_{k=0}^{n}\sum_{l=0}^{k}\mathcal{E}_{l,\lambda}^{(r)}(x)S_{J,\lambda}^{(1)}(k,l)S_{2,\lambda}(n,k).
    \end{displaymath}
\end{corollary}
In \eqref{16}, replacing $t$ by $e_{\lambda}(e_{\lambda}(t)-1)-1$, we get
\begin{align}
&\bigg(\frac{2}{e_{\lambda}\big(e_{\lambda}(t)-1\big)+e_{\lambda}^{-1}\big(e_{\lambda}(t)-1\big)}\bigg)^{r}e_{\lambda}^{x}\big(e_{\lambda}(t)-1\big)\label{25} \\
&=\ \sum_{k=0}^{\infty}C_{k}^{(r)}(x)\frac{1}{k!}\big(e_{\lambda}(e_{\lambda}(t)-1)-1\big)^{k}\nonumber \\
&=\ \sum_{k=0}^{\infty}C_{k}^{(r)}(x)\sum_{n=k}^{\infty}S_{J,\lambda}^{(2)}(n,k)\frac{t^{n}}{n!}\nonumber \\
&=\ \sum_{n=0}^{\infty}\bigg(\sum_{k=0}^{n}C_{k}^{(r)}(x)S_{J,\lambda}^{(2)}(n,k)\bigg)\frac{t^{n}}{n!}.\nonumber
\end{align}

On the other hand,
\begin{align}
&\bigg(\frac{2}{e_{\lambda}\big(e_{\lambda}(t)-1\big)+e_{\lambda}^{-1}\big(e_{\lambda}(t)-1\big)}\bigg)^{r}e_{\lambda}^{x}\big(e_{\lambda}(t)-1\big) \label{26} \\
&=\sum_{k=0}^{\infty}\mathcal{E}_{k,\lambda}^{(r)}(x)\frac{1}{k!}\big(e_{\lambda}(t)-1\big)^{k} \nonumber \\
&= \sum_{k=0}^{\infty}\mathcal{E}_{k,\lambda}^{(r)}(x)\sum_{n=k}^{\infty}S_{2,\lambda}(n,k)\frac{t^{n}}{n!}\nonumber \\
&= \sum_{n=0}^{\infty}\bigg(\sum_{k=0}^{n}\mathcal{E}_{k,\lambda}^{(r)}(x)S_{2,\lambda}(n,k)\bigg)\frac{t^{n}}{n!}. \nonumber
\end{align}

Therefore, by \eqref{25} and \eqref{26}, we obtain the following theorem.
\begin{theorem}
    For $n\ge 0$, we have
    \begin{displaymath}
        \sum_{k=0}^{n}S_{J,\lambda}^{(2)}(n,k)C_{k}^{(r)}(x)\ =\ \sum_{k=0}^{n} \mathcal{E}_{k,\lambda}^{(r)}(x)S_{2,\lambda}(n,k).
    \end{displaymath}
\end{theorem}
From Theorem 2, we have the following corollary.
\begin{corollary}
    For $n\ge 0$, we have
    \begin{displaymath}
        \mathcal{E}_{n,\lambda}^{(r)}(x)\ =\ \sum_{k=0}^{n}\sum_{l=0}^{k}S_{1,\lambda}(n,k)S_{J,\lambda}^{(2)}(k,l)C_{l}^{(r)}(x).
    \end{displaymath}
\end{corollary}
From \eqref{9}, we note that
\begin{align}
&\bigg(\frac{2}{e_{\lambda}\big(e_{\lambda}(t)-1\big)+e_{\lambda}^{-1}\big(e_{\lambda}(t)-1\big)}\bigg)^{r}e_{\lambda}^{x}\big(e_{\lambda}(t)-1\big) \label{27} \\
&=\ \sum_{l=0}^{\infty}\mathcal{E}_{l,\lambda}^{(r)}\frac{1}{l!}\big(e_{\lambda}(t)-1\big)^{l}\sum_{m=0}^{\infty}B_{m,\lambda}(x)\frac{t^{m}}{m!}\nonumber \\
&=\ \sum_{k=0}^{\infty}\sum_{l=0}^{k}\mathcal{E}_{l,\lambda}^{(r)}S_{2,\lambda}(k,l)\frac{t^{k}}{k!}\sum_{m=0}^{\infty}B_{m,\lambda}(x)\frac{t^{m}}{m!} \nonumber\\
&=\ \sum_{n=0}^{\infty}\bigg(\sum_{k=0}^{n}\binom{n}{k}\sum_{l=0}^{k}\mathcal{E}_{l,\lambda}^{(r)}S_{2,\lambda}(k,l)B_{n-k,\lambda}(x)\bigg)\frac{t^{n}}{n!}. \nonumber
\end{align}
Therefore, by \eqref{25} and \eqref{27}, we obtain the following theorem.
\begin{theorem}
For $n\ge 0$, we have
\begin{displaymath}
\sum_{k=0}^{n}S_{J,\lambda}^{(2)}(n,k)C_{k}^{(r)}(x)\ =\ \sum_{k=0}^{n}\binom{n}{k}\sum_{l=0}^{k}\mathcal{E}_{l,\lambda}^{(r)}S_{2,\lambda}(k,l)B_{n-k,\lambda}(x).
\end{displaymath}
\end{theorem}

\section{Discrete harmonic numbers and related polynomials and numbers}

From \eqref{6}, we note that
\begin{equation}
\lim_{\lambda\rightarrow 0}\log_{\lambda}(1+t)\ =\ \sum_{n=1}^{\infty}\frac{(-1)^{n-1}}{n}t^{n}=\log(1+t). \label{29}
\end{equation}
In view of \eqref{8} and \eqref{29}, we may consider the degenerate harmonic numbers given by
\begin{equation}
-\frac{\log_{\lambda}(1-t)}{1-t}\ =\ \sum_{n=1}^{\infty}H_{n,\lambda}t^n. \label{30}
\end{equation}
Note that
\begin{displaymath}
\lim_{\lambda\rightarrow 0}H_{n,\lambda}\ =\ H_{n}\ =\ 1+\frac{1}{2}+\cdots+\frac{1}{n},\quad(n\in\mathbb{N}).
\end{displaymath}
From \eqref{6} and \eqref{30}, we note that
\begin{align}
\frac{-\log_{\lambda}(1-t)}{1-t}\
&=\ \sum_{m=0}^{\infty}t^{m}\sum_{k=1}^{\infty}(-\lambda)^{k-1}(1)_{k,1/\lambda}\frac{t^{k}}{k!} \label{31} \\
&=\ \sum_{n=1}^{\infty}\bigg(\sum_{k=1}^{n}\frac{(-\lambda)^{k-1}(1)_{k,1/\lambda}}{k!}\bigg)t^{n}. \nonumber
\end{align}
Thus, by \eqref{30} and \eqref{31}, we get
\begin{displaymath}
H_{n,\lambda}=\sum_{k=1}^{n}\frac{(-\lambda)^{k-1}(1)_{k,1/\lambda}}{k!},\quad(n\in\mathbb{N}).
\end{displaymath}
Indeed,
\begin{displaymath}
\lim_{\lambda\rightarrow 0}H_{n,\lambda}\ =\ 1+\frac{1}{2}+\frac{1}{3}+\cdots+\frac{1}{n}=H_{n}.
\end{displaymath}
For $r\in\mathbb{N}$, the degenerate harmonic polynomials $H_{n,\lambda}^{(r)}(x)$ of order $r$ are defined by
\begin{equation}
\sum_{n=0}^{\infty}H_{n,\lambda}^{(r)}(x)t^{n}\ =\ \frac{\big(-\log_{\lambda}(1-t)\big)^{r+1}}{t(1-t)}(1-t)^{x}. \label{32}
\end{equation}
When $x=0$, $H_{n,\lambda}^{(r)}=H_{n,\lambda}^{(r)}(0)$ are called the degenerate harmonic numbers of order $r$. \\
From \eqref{32}, we note that
\begin{equation}\label{33}
\sum_{n=0}^{\infty}H_{n,\lambda}^{(0)}t^{n}\ =\ \frac{-\log_{\lambda}(1-t)}{t(1-t)}
=\frac{1}{t}\sum_{n=1}^{\infty}H_{n,\lambda}t^n
=\ \sum_{n=0}^{\infty}H_{n+1,\lambda}t^{n}.
\end{equation}
Comparing the coefficients on both sides of \eqref{33}, we obtain
\begin{equation}
H_{n,\lambda}^{(0)}=H_{n+1,\lambda},\quad(n\ge 0). \label{34}
\end{equation}
From \eqref{32} and \eqref{34}, we have
\begin{align}
\sum_{n=0}^{\infty}H_{n,\lambda}^{(0)}(x)t^{n}&=\frac{-\log_{\lambda}(1-t)}{t(1-t)}(1-t)^{x}\label{35}\\
&=\sum_{m=0}^{\infty}H_{m+1,\lambda}t^{m}\sum_{l=0}^{\infty}\binom{x}{l}(-1)^{l}x^{l} \nonumber \\
&= \sum_{n=0}^{\infty}\bigg(\sum_{m=0}^{n}H_{m+1,\lambda}\binom{x}{n-m}(-1)^{n-m}\bigg)t^{n}. \nonumber
\end{align}
By \eqref{35}, we get
\begin{displaymath}
    H_{n,\lambda}^{(0)}(x)\ =\ \sum_{m=0}^{n}H_{m+1,\lambda}\binom{x}{n-m}(-1)^{n-m},\quad (n\ge 0).
\end{displaymath}
Now, we observe that
\begin{align}
&\frac{\big(-\log_{\lambda}(1-t)\big)^{r+1}}{t(1-t)}=\frac{1}{t(1-t)}\bigg(\sum_{l=1}^{\infty}\frac{(-\lambda)^{l-1}}{l!}(1)_{l,1/\lambda}t^{l}\bigg)^{r+1}\label{36} \\
&=\ \bigg(\frac{1}{t}+\frac{1}{1-t}\bigg)\sum_{l=r+1}^{\infty}\sum_{l_{1}+\cdots+l_{r+1}=l}\frac{(-\lambda)^{l-r-1}}{l_{1}!l_{2}!\cdots l_{r+1}!} (1)_{l_{1},1/\lambda}\cdots(1)_{l_{r+1},1/\lambda}\,t^{l} \nonumber \\
&=\ \sum_{n=r}^{\infty} \sum_{l_{1}+\cdots+l_{r+1}=n+1}\frac{(-\lambda)^{n-r}(1)_{l_{1},1/\lambda}\cdots(1)_{l_{r+1},1/\lambda}}{l_{1}!\cdots l_{r+1}!}\,t^{n} \nonumber\\
&\quad +\sum_{n=r+1}^{\infty}\sum_{l=r+1}^{n}\sum_{l_{1}+\cdots+l_{r+1}=l}\frac{(-\lambda)^{l-r-1}(1)_{l_{1},1/\lambda}\cdots(1)_{l_{r+1},1/\lambda}}{l_{1}!\cdots l_{r+1}!}\,t^{n}\nonumber \\
&=\sum_{n=r}^{\infty}\bigg(\sum_{l=r+1}^{n+1}\sum_{l_{1}+\cdots+l_{r+1}=l}\frac{(-\lambda)^{l-r-1}}{l_{1}!l_{2}!\cdots l_{r+1}!}(1)_{l_{1},1/\lambda}(1)_{l_{2},1/\lambda}\cdots (1)_{l_{r+1},1/\lambda}\bigg)t^{n}\nonumber.
\end{align}
Therefore, by \eqref{32} and \eqref{36}, we obtain the following equation.
\begin{align}
& H_{n,\lambda}^{(r)}=\left\{\begin{array}{ccc}
\displaystyle\sum_{l=r+1}^{n+1}\sum_{l_{1}+\cdots+l_{r+1}=l}\frac{(-\lambda)^{l-r-1}}{l_{1}!l_{2}!\cdots l_{r+1}!}(1)_{l_{1},1/\lambda}(1)_{l_{2},1/\lambda}, \cdots(1)_{l_{r+1},1/\lambda},\displaystyle & \textrm{if $n\ge r$}, \\
0, & \textrm{otherwise}.
\end{array}\right.\label{37}
\end{align}
By \eqref{32} and \eqref{36}, we get
\begin{align}
\sum_{n=0}^{\infty}H_{n,\lambda}^{(r)}(x)t^{n}&= \frac{\big(-\log_{\lambda}(1-t)\big)^{r+1}}{t(1-t)}(1-t)^{x} \label{38} \\
&=\sum_{k=r}^{\infty}\bigg(\sum_{l=r+1}^{k+1}\sum_{l_{1}+\cdots+l_{r+1}=l}\frac{(-\lambda)^{l-r-1}(1)_{l_{1},1/\lambda}\cdots (1)_{l_{r+1},1/\lambda}}{l_{1}!\cdots l_{r+1}!}\bigg)t^{k}\sum_{m=0}^{\infty}\binom{x}{m}(-1)^{m}t^{m}\nonumber \\
&=\sum_{n=r}^{\infty}\bigg(\sum_{k=r}^{n}\binom{x}{n-k}(-1)^{n-k}\sum_{l=r+1}^{k+1}\sum_{l_{1}+\cdots+l_{r+1}=l}\frac{(-\lambda)^{l-r-1}(1)_{l_{1},1/\lambda}\cdots (1)_{l_{r+1},1/\lambda}}{l_{1}!\cdots l_{r+1}!}\bigg)t^{n}.\nonumber
\end{align}
By comparing the coefficients on both sides of \eqref{38}, we get
\begin{equation}
H_{n,\lambda}^{(r)}(x)=\left\{\begin{array}{cc}
\displaystyle\sum_{k=r}^{n}\binom{x}{n-k}(-1)^{n-k}\sum_{l=r+1}^{k+1}\sum_{l_{1}+\cdots+l_{r+1}=l}\frac{(-\lambda)^{l-r-1}(1)_{l_{1},1/\lambda}\cdots (1)_{l_{r+1},1/\lambda}}{l_{1}!\cdots l_{r+1}!}\displaystyle, & \textrm{if $n\ge r$,}\\
0 , & \textrm{otherwise.}
\end{array}\right. \label{39}
\end{equation}
Let us consider the generalized degenerate harmonic numbers which are given by
\begin{equation}
    \sum_{n=0}^{\infty}H_{\lambda}(n+r+1,r)t^{n}\ =\ \frac{\big(-\log_{\lambda}(1-t)\big)^{r+1}}{t^{r+1}(1-t)}. \label{40}
\end{equation}
Replacing $t$ by $-t$, we get
\begin{equation}
    \sum_{n=0}^{\infty}H_{\lambda}(n+r+1,r)(-t)^{n}\ =\ \frac{\big(\log_{\lambda}(1+t)\big)^{r+1}}{t^{r+1}(1+t)}.\label{41}
\end{equation}
From \eqref{6}, we note that
\begin{align}
\big(\log_{\lambda}(1+t)\big)^{r+1}\ &=\ \bigg(\sum_{l=1}^{\infty}\lambda^{l-1}(1)_{l,1/\lambda}\frac{t^{l}}{l!}\bigg)^{r+1}\label{42} \\
&=\ \sum_{l=r+1}^{\infty}\bigg(\sum_{l_{1}+\cdots+l_{r+1}=l}\frac{\lambda^{l-r-1}(1)_{l_{1},1/\lambda}\cdots(1)_{l_{r+1},1/\lambda}}{l_{1}!l_{2}!\cdots l_{r+1}!}
\bigg)t^{l}.\nonumber
\end{align}
By \eqref{42}, we get
\begin{align}
\frac{\big(\log_{\lambda}(1+t)\big)^{r+1}}{t^{r+1}(1+t)}&=\frac{1}{t^{r+1}}\sum_{l=r+1}^{\infty}\sum_{l_{1}+\cdots+l_{r+1}=l}\frac{\lambda^{l-r-1}(1)_{l_{1},1/\lambda}\cdots(1)_{l_{r+1},1/\lambda}}{l_{1}!\cdots l_{r+1}!}t^{l}\sum_{m=0}^{\infty}(-1)^{m}t^{m}\label{43} \\
&=\frac{1}{t^{r+1}}\sum_{n=r+1}^{\infty}\sum_{l=r+1}^{n}(-1)^{n-l}\sum_{l_{1}+\cdots+l_{r+1}=l}\frac{\lambda^{l-r-1}(1)_{l_{1},1/\lambda}\cdots(1)_{l_{r+1},1/\lambda}}{l_{1}!\cdots l_{r+1}!}\,t^{n}\nonumber \\
&=\sum_{n=0}^{\infty}\bigg(\sum_{l=r+1}^{n+r+1}(-1)^{n+r+1-l}\sum_{l_{1}+\cdots+l_{r+1}=l}\frac{\lambda^{l-r-1}(1)_{l_{1},1/\lambda}\cdots(1)_{l_{r+1},1/\lambda}}{l_{1}!\cdots l_{r+1}!}\bigg)t^{n}.\nonumber
\end{align}
Therefore, by \eqref{41} and \eqref{43}, we get
\begin{equation}
    H_{\lambda}(n+r+1,r)=\sum_{l=r+1}^{n+r+1}(-1)^{r+1-l} \sum_{l_{1}+\cdots+l_{r+1}=l}\frac{\lambda^{l-r-1}(1)_{l_{1},1/\lambda}\cdots(1)_{l_{r+1},1/\lambda}}{l_{1}!\cdots l_{r+1}!}.\label{44}
\end{equation}

Let $r$ be a positive integer. We observe that
\begin{align}
&\frac{\big(\log_{\lambda}(1+t)\big)^{r+1}}{t^{r+1}(1+t)}\bigg(\frac{2}{(1+t)+(1+t)^{-1}}\bigg)^{r}(1+t)^{x}\label{45}\\
&\quad=\bigg(\frac{2}{(1+t)+(1+t)^{-1}}\bigg)^{r}(1+t)^{x-1}\frac{1}{t^{r+1}}\big(\log_{\lambda}(1+t)\big)^{r+1}\nonumber\\
&\quad=\sum_{l=0}^{\infty}C_{l}^{(r)}(x-1)\frac{t^{l}}{l!}\sum_{k=0}^{\infty}S_{1,\lambda}(k+r+1,r+1)\frac{k!(r+1)!}{(k+r+1)!}\frac{t^{k}}{k!} \nonumber \\
&\quad=\sum_{l=0}^{\infty}C_{l}^{(r)}(x-1)\frac{t^{l}}{l!}\sum_{k=0}^{\infty}\frac{S_{1,\lambda}(k+r+1,r+1)}{\binom{k+r+1}{k}}\frac{t^{k}}{k!} \nonumber\\
&\quad=\sum_{n=0}^{\infty}\bigg(\sum_{k=0}^{n}\frac{\binom{n}{k}S_{1,\lambda}(k+r+1,r+1)}{\binom{k+r+1}{k}}C_{n-k}^{(r)}(x-1)\bigg)\frac{t^{n}}{n!}.\nonumber
\end{align}

On the other hand,
\begin{align}
&\frac{\big(\log_{\lambda}(1+t)\big)^{r+1}}{t^{r+1}(1+t)}\bigg(\frac{2}{(1+t)+(1+t)^{-1}}\bigg)^{r}(1+t)^{x} \label{46} \\
&\quad=\sum_{l=0}^{\infty}H_{\lambda}(l+r+1)(-t)^{l}\sum_{m=0}^{\infty}C_{m}^{(r)}(x)\frac{t^{m}}{m!}\nonumber \\
&\quad=\sum_{n=0}^{\infty}\bigg(\sum_{l=0}^{n}\binom{n}{l}l!H_{\lambda}(l+r+1,r)(-1)^{l}C_{n-l}^{(r)}(x)\bigg)\frac{t^{n}}{n!}.\nonumber
\end{align}
Therefore, by \eqref{45} and \eqref{46}, we obtain the following equation.
\begin{equation}
    \sum_{k=0}^{n}\frac{\binom{n}{k}S_{1,\lambda}(k+r+1,r+1)}{\binom{k+r+1}{k}}C_{n-k}^{(r)}(x-1)=\sum_{l=0}^{n}\binom{n}{l}l!H_{\lambda}(l+r+1,r)(-1)^{l}C_{n-l}^{(r)}(x).\label{47}
\end{equation}

\section{Conclusion}

In Section 2, the orthogonality relations were shown for the degenerate Stirling numbers from
which the inversion theorem was derived. Then, in connection with the degenerate Stirling
numbers of both kinds and the Jindalrae-Stirling numbers of both kinds, several relations
were proved between the type 2 degenerate Euler polynomials of order $r$ and the type 2
Changhee polynomials of order $r$. In additon, a recurrence relation was deduced for the
type 2 higher-order degenerate Euler polynomials. In Section 3, the discrete harmonic
numbers were introduced as a degenerate version of the usual harmonic numbers. Then,
as a natural generalization of these numbers, the higher-order degenerate harmonic polynomials
were considered and an explicit expression for them was obtained. Finally, the generalized
degenerate harmonic numbers were constructed so that an identity involving these numbers,
the type 2 higher-order Changhee polynomials and the degenerate Stirling numbers of the first kind was derived. \\
\indent We would like to mention three possibilities for applications of our results to other areas.
 The first possibility is their applications to differential equations. In [8], certain infinite families of ordinary
 differential equations, satisfied by the generating functions of some degenerate polynomials, were derived in order
 to find new combinatorial identities for those polynomials. The second possibility is their applications to
 probability theory. In [15,18], by using the generating functions of the moments of certain random variables,
 new identities connecting some special numbers and moments of random variables were deduced.
The third possibility is their applications to identities of symmetry. In [14], abundant identities of
symmetry were obtained  for various degenerate versions of many special polynomials by using
$p$-adic fermionic integrals. \\
\indent It is one of our future projects to continue to study various degenerate versions of some special
polynomials and numbers, and to find some of their possible applications to mathematics, science and engineering.

\end{document}